\documentclass{amsart}

\usepackage{amssymb}

\allowdisplaybreaks

\newtheorem{theorem}{Theorem}[section]

\newcommand{\dem}{\noindent {\bf Proof. }}
\newcommand{\fin}{\hspace*{\fill} $\square$ \vskip0.2cm}

\newcommand{\C}{\mathbb{C}}

\newcommand{\G}{\widehat{\mathrm{G}}}

\begin{document}

\title[Almost unimodular systems] {Almost unimodular systems \\ on
compact groups with disjoint spectra}

\author[Javier Parcet]
{Javier Parcet}

\address{Universidad Aut\'{o}noma de Madrid}

\email{javier.parcet@uam.es}

\footnote{Partially supported by the Project BFM 2001/0189,
Spain.}

\date{}

\begin{abstract}
We emulate the Rademacher functions on any non-commutative compact
group requiring the resulting system to have pairwise disjoint
spectra.
\end{abstract}

\maketitle

\section*{Introduction}

Let $\mathrm{G}$ be a compact topological group with normalized
Haar measure $\mu$ and dual object $\G$. Let $\pi: \mathrm{G}
\rightarrow U(\mathcal{H}_{\pi})$ be an irreducible unitary
representation of $\mathrm{G}$. Then, given an integrable function
$f: \mathrm{G} \to \C$, the Fourier coefficient of $f$ at $ \pi$
is defined as follows $$\widehat{f}(\pi) = \int_{\mathrm{G}} f(g)
\pi(g)^* \, d \mu(g) \in \mathcal{B}(\mathcal{H}_{\pi}).$$ It is
well-known how to construct $n$ functions $f_1, f_2, \ldots f_n:
\mathrm{G} \to \C$ satisfying
\begin{itemize}
\item[$(\mathrm{P}1)$] $f_1,
f_2, \ldots, f_n$ have pairwise disjoint supports on $\mathrm{G}$.
\item[$(\mathrm{P}2)$] The norm of $\widehat{f}_k(\pi)$ on
$\mathcal{S}^{d_{\pi}}_{q'}$ does not depend on $k$ for any $\pi
\in \G$.
\end{itemize}
Just take disjoint translations of a common function on
$\mathrm{G}$ with sufficiently small support. Also, one can
consider the dual properties with respect to the Fourier transform
on $\mathrm{G}$. Namely,
\begin{itemize}
\item[$(\widehat{\mathrm{P}1})$] $\widehat{f}_1,
\widehat{f}_2, \ldots, \widehat{f}_n$ have pairwise disjoint
supports on $\G$.
\item[$(\widehat{\mathrm{P}2})$] The absolute value $|f_k(g)|$
does not depend on $k$ for almost all $g \in \mathrm{G}$.
\end{itemize}

\vskip3pt

If $\mathrm{G}$ is abelian we can easily construct such a system
by taking $n$ irreducible characters of $\mathrm{G}$. Moreover,
many other constructions are available. Namely, we can take $n$
disjoint translations of a common function on $\G$. Then we get
the desired property by Pontrjagin duality. This kind of systems
have been applied to study the Fourier type constants of
finite-dimensional Lebesgue spaces: if $1 \le p < q \le 2$, the
Fourier $q$-type constant of $\ell_p(n)$ with respect to $\G$ is
$n^{1/p - 1/q}$, which is optimal among the family of $n$
dimensional Banach spaces. In particular, this provides sharp
results about the Fourier type of infinite-dimensional Lebesgue
spaces.

\vskip3pt

However, irreducible characters on non-commutative compact groups
are no longer unimodular. Furthermore, the dual object is not a
group anymore and consequently Pontrjagin duality does not hold.
Besides, Tannaka's theorem (the non-commutative counterpart of
Pontrjagin theorem) does not fit properly in this context. Hence,
it is natural to wonder if it is possible to construct a system
$$\Phi = \Big\{ f_m: \mathrm{G} \to \C \, \big| \ m \ge 1 \Big\}$$
made up of functions $f_1, f_2, \ldots$ satisfying similar
properties to $(\widehat{\mathrm{P}1})$ and
$(\widehat{\mathrm{P}2})$. In this paper, we construct an
\emph{almost} unimodular system of trigonometric polynomials $f_m:
\mathrm{G} \to \C$ with non-overlapping ranges of frequencies. The
work \cite{GMP} is the main motivation for this note.

\section{Almost unimodular systems}

\begin{theorem} \label{Non-overlapping}
Let $\mathrm{G}$ be an infinite compact group and let $f_0$ be a
continuous function in $L_2(\mathrm{G})$. Let $\varepsilon_1,
\varepsilon_2, \ldots$ be any sequence of positive numbers
decreasing to $0$. Then there exists a collection $\Omega_1,
\Omega_2, \ldots$ of measurable subsets of $\mathrm{G}$ such that
$\mu(\Omega_m) \rightarrow 1$ as $m \rightarrow \infty$, and a
system $\Phi = \big\{ f_m \, | \ m \ge 1 \big\}$ of trigonometric
polynomials in $\mathrm{G}$ satisfying:
\begin{itemize}
\item[\textnormal{i)}] $\widehat{f}_1, \widehat{f}_2,
\ldots$ have pairwise disjoint supports on $\G$.
\item[\textnormal{ii)}] $|f_m(g)| < |f_0(g)| + \varepsilon_m$
for all $g \in \mathrm{G}$ and all $m \ge 1$.
\item[\textnormal{iii)}] $|f_m(g)| > |f_0(g)| - \varepsilon_m$
for all $g \in \Omega_m$ and all $m \ge 1$.
\end{itemize}
If $f_0 \equiv 1$ we have an almost unimodular system on
$\mathrm{G}$ with pairwise disjoint spectra.
\end{theorem}

\dem First we recall that it is essentially no restriction to
assume that a compact topological group is Hausdorff, see e.g.
Corollary 2.3 in \cite{Fo}. Then we point out that the normalized
Haar measure $\mu$ of an infinite Hausdorff compact group
$\mathrm{G}$ has no atoms. This is an easy consequence of the
translation invariance and finiteness of $\mu$. Indeed, let us
assume that $\mu$ has an atom $\Omega$ with $\mu(\Omega) = \alpha$
so that $0 < \alpha \le 1$. Let us consider an open neighborhood
$\mathcal{U}$ of $\mathrm{G}$. Since $\mathrm{G}$ is non-finite
and Hausdorff, we can choose $\mathcal{U}$ small enough so that
there exists $m_0$ pairwise disjoint translations $$g_1
\mathcal{U}, g_2 \mathcal{U}, \ldots, g_{m_0} \mathcal{U}$$ with
$m_0 > 1/\alpha$. On the other hand, the compactness of
$\mathrm{G}$ allows us to write $$\mathrm{G} = \bigcup_{k=1}^{n_0}
h_k \mathcal{U}$$ as a finite union of translations of
$\mathcal{U}$. In particular, there must exists $1 \le k_0 \le
n_0$ such that $\mu(h_{k_0} \mathcal{U} \cap \Omega) > 0$.
Moreover, we must have $\mu(h_{k_0} \mathcal{U} \cap \Omega) =
\alpha$ since $\Omega$ is an atom. Then, by the translation
invariance of the Haar measure, we have $$\mu(g_j \mathcal{U}) \ge
\mu(g_j \mathcal{U} \cap g_j h_{k_0}^{-1} \Omega) = \mu(h_{k_0}
\mathcal{U} \cap \Omega) = \alpha.$$ Then, since $g_1 \mathcal{U},
g_2 \mathcal{U}, \ldots, g_{m_0} \mathcal{U}$ are pairwise
disjoint, we obtain $$\mu(\mathrm{G}) \ge \sum_{j=1}^{m_0} \mu(g_j
\mathcal{U}) \ge m_0 \alpha > 1,$$ which contradicts the
assumption that $\mu$ has mass $1$ and proves our claim.

\vskip3pt

\noindent \textbf{Step 1.} Since the Haar measure $\mu$ has no
atoms, it is clear that the same holds for the measure $d \nu =
|f_0|^2 d \mu$. Then, by the absence of $\nu$-atoms we can define
a family of $\nu$-measurable dyadic sets $$\Big\{ \mathrm{D}_j^k
\, \big| \ k \ge 1, \, 1 \le j \le 2^k \Big\}$$ in $\mathrm{G}$
satisfying the standard conditions:
\begin{itemize}
\item[(a)] $\displaystyle \mathrm{D}_j^k = \mathrm{D}_{2j-1}^{k+1}
\cup \mathrm{D}_{2j}^{k+1}$ for all $j,k$.
\item[(b)] $\mathrm{D}_1^k, \mathrm{D}_2^k, \ldots, \mathrm{D}_{2^k}^k$
are pairwise disjoint and $$\mathrm{G} = \bigcup_{j=1}^{2^k}
\mathrm{D}_j^k \quad \mbox{for all} \quad k \ge 1.$$
\item[(c)] The sets $\mathrm{D}_j^k$ have the same $\nu$-measure for
$k$ fixed $$\int_{\mathrm{D}_j^k} |f_0(g)|^2 d\mu(g) = 2^{-k} \,
\int_{\mathrm{G}} |f_0(g)|^2 d\mu(g).$$
\end{itemize}
Then, if $1_{\Omega}$ denotes the characteristic function of
$\Omega$, we define the system $$\Delta = \Big\{ \delta_k:
\mathrm{G} \to \C \, \big| \ k \ge 1 \Big\}$$ as follows
$$\delta_k(g) = f_0(g) \, \sum_{j=1}^{2^k} (-1)^{j+1}
1_{\mathrm{D}_j^k}(g).$$ Since every finite Radon measure is
regular, we can consider compact sets $\mathrm{K}_j^k \subset
\mathrm{D}_j^k$ such that $\mu(\mathrm{D}_j^k \setminus
\mathrm{K}_j^k) < 2^{-2k}$ for all $k \ge 1$ and all $1 \le j \le
2^k$. Then, by Urylshon's lemma we can construct continuous
functions $\gamma_j^k$ on $\mathrm{G}$ such that
$1_{\mathrm{K}_j^k} \le \gamma_j^k \le 1_{\mathrm{D}_j^k}$. Then,
we define the system $$\Psi = \Big\{ \psi_k: \mathrm{G} \to \C \,
\big| \ k \ge 1 \Big\}$$ as follows $$\psi_k(g) = f_0(g) \,
\sum_{j=1}^{2^k} (-1)^{j+1} \gamma_j^k(g).$$

\vskip3pt

\noindent \textbf{Step 2.} Let $\pi: \mathrm{G} \rightarrow
\mathcal{B}(\mathcal{H}_{\pi})$ be an irreducible unitary
representation of degree $d_{\pi}$ and let us fix $1 \le i,j \le
d_{\pi}$. Since $\pi(g)$ can be identified with a unitary $d_{\pi}
\times d_{\pi}$ matrix, we have
\begin{eqnarray*}
\lefteqn{\Big( \sum_{k=1}^{\infty} |\widehat{\psi}_k(\pi)_{ij}|^2
\Big)^{1/2}} \\ & = & \Big( \sum_{k=1}^{\infty} \Big|
\int_{\mathrm{G}} \psi_k(g) \overline{\pi_{ji}(g)} \, d \mu(g)
\Big|^2 \Big)^{1/2} \\ & \le & \Big( \sum_{k=1}^{\infty} \Big|
\int_{\mathrm{G}} (\psi_k - \delta_k)(g) \overline{\pi_{ji}(g)} d
\mu(g) \Big|^2 \Big)^{1/2} + \Big( \sum_{k=1}^{\infty} \Big|
\int_{\mathrm{G}} \delta_k(g) \overline{\pi_{ji}(g)} d \mu(g)
\Big|^2 \Big)^{1/2}
\end{eqnarray*}
Applying H\"{o}lder inequality for the first part and Bessel
inequality for the second (note that the system
$\|f_0\|_{L_2(\mathrm{G})}^{-1} \Delta$ is orthonormal in
$L_2(\mathrm{G})$) we get
\begin{eqnarray*}
\Big( \sum_{k=1}^{\infty} |\widehat{\psi}_k(\pi)_{ij}|^2
\Big)^{1/2} \!\! & \le & \!\! \|\pi_{ji}\|_{L_2(\mathrm{G})} \Big(
\sum_{k=1}^{\infty} \|\psi_k - \delta_k\|_{L_2(\mathrm{G})}^2
\Big)^{1/2} + \|f_0\|_{L_2(\mathrm{G})}
\|\pi_{ji}\|_{L_2(\mathrm{G})}
\\ \!\! & \le & \!\! \frac{1}{\sqrt{d_{\pi}}} \, \Big(
\|f_0\|_{L_{\infty}(\mathrm{G})} + \|f_0\|_{L_2(\mathrm{G})} \Big)
< \infty.
\end{eqnarray*}
The last inequality uses the estimate $$\|\psi_k -
\delta_k\|_{L_2(\mathrm{G})}^2 \le \sum_{j=1}^{2^k}
\int_{\mathrm{D}_j^k \setminus \mathrm{K}_j^k} |f_0(g)|^2 d \mu(g)
\le 2^{-k} \|f_0\|_{L_{\infty}(\mathrm{G})}^2,$$ which follows
since $|\psi_k - \delta_k|^2$ is supported in
$$\bigcup_{j=1}^{2^k} \mathrm{D}_j^k \setminus \mathrm{K}_j^k$$
and bounded above by $|f_0|^2$. In particular,
$|\widehat{\psi}_k(\pi)_{ij}| \to 0$ as $ k \to \infty$.
Therefore, given any $\delta > 0$ and any finite subset $\Lambda
\subset \G$, there exists a positive integer
$\mathrm{M}(\Lambda,\delta)$ such that for all $k \ge
\mathrm{M}(\Lambda,\delta)$ we have
\begin{equation} \label{delta}
\sum_{\pi \in \Lambda} d_{\pi} \sum_{i,j = 1}^{d_{\pi}}
|\widehat{\psi}_k(\pi)_{ij}| < \delta.
\end{equation}

\noindent \textbf{Step 3.} For any $\pi \in \G$ let us denote by
$\mathcal{E}_{\pi}$ the linear span of the entries of $\pi$. Also,
we shall write $\mathcal{E}$ for the linear span of the union of
the spaces $\mathcal{E}_{\pi}$ when $\pi$ runs over $\G$. That is,
$\mathcal{E}$ is the space of trigonometric polynomials on
$\mathrm{G}$. By the Peter-Weyl theorem we know that $\mathcal{E}$
is dense in the space $\mathcal{C}(\mathrm{G})$ of continuous
functions on $\mathrm{G}$ with respect to the uniform norm. Then
we construct the system $$\Phi = \Big\{ f_m: \mathrm{G} \to \C \,
\big| \ m \ge 1 \Big\}$$ as follows:
\begin{itemize}
\item[1.] Let $f_1 \in \mathcal{E}$ be such that
$\displaystyle \sup_{g \in \mathrm{G}} |f_1(g) - \psi_1(g)| <
\varepsilon_1$.
\item[2.] For $m > 1$ we proceed by induction. Let $$\Lambda_m =
\bigcup_{k=1}^{m-1} \mbox{supp} \, \widehat{f}_k \subset \G \quad
\mbox{and} \quad \delta_m = \min \Big( \frac{\varepsilon_m}{3},
\Big[ \sum_{\pi \in \Lambda_m} d_{\pi}^{5/2} \Big]^{-1} \Big).$$
Let $k_m = \mathrm{M}(\Lambda_m, \delta_m)$ and let $\xi_m \in
\mathcal{E}$ be such that $$\sup_{g \in \mathrm{G}} |\xi_m(g) -
\psi_{k_m}(g)| < \delta_m^2.$$ Then we define $$f_m(g) = \xi_m(g)
- \sum_{\pi \in \Lambda_m} d_{\pi} \mbox{tr}(\pi(g)
\widehat{\xi}_m(\pi)).$$
\end{itemize}

\vskip3pt

\noindent \textbf{Step 4.} Kunze's Hausdorff-Young inequality and
(\ref{delta}) give
\begin{eqnarray*}
\sup_{g \in \mathrm{G}} |f_m(g) - \psi_{k_m}(g)| & \le & \sup_{g
\in \mathrm{G}} |\xi_m(g) - \psi_{k_m}(g)| + \sup_{g \in
\mathrm{G}} \Big| \sum_{\pi \in \Lambda_m} d_{\pi}
\mbox{tr}(\pi(g) \widehat{\xi}_m(\pi)) \Big| \\ & < & \delta_m^2 +
\sum_{\pi \in \Lambda_m}^{\null} d_{\pi}
\|\widehat{\xi}_m(\pi)\|_{\mathcal{S}^{d_{\pi}}_1} \\ & \le &
\delta_m^2 + \sum_{\pi \in \Lambda_m} d_{\pi}
\sum_{i,j=1}^{d_{\pi}} |\widehat{\xi}_m(\pi)_{ij}| \\ & < &
\delta_m^2 + \delta_m + \sum_{\pi \in \Lambda_m} d_{\pi}
\sum_{i,j=1}^{d_{\pi}} \Big| \int_{\mathrm{G}} (\xi_m -
\psi_{k_m})(g) \overline{\pi_{ji}(g)} \, d \mu(g) \Big| \\ & \le &
\delta_m^2 + \delta_m + \|\xi_m - \psi_{k_m}\|_{L_2(\mathrm{G})}
\sum_{\pi \in \Lambda_m} d_{\pi} \sum_{i,j=1}^{d_{\pi}}
\|\pi_{ij}\|_{L_2(\mathrm{G})} \\ & \le & \delta_m^2 + \delta_m +
\sup_{g \in \mathrm{G}} |\xi_m(g) - \psi_{k_m}(g)| \sum_{\pi \in
\Lambda_m}^{\null} d_{\pi}^{5/2} \\ & < & \delta_m^2 + 2 \delta_m
\le \varepsilon_m.
\end{eqnarray*}
We finally show that $\Phi$ satisfies all the properties stated
above. Taking $$\Omega_m = \bigcup_{1 \le j \le 2^{k_m}}
\mathrm{K}_j^{k_m},$$ we have $$\mu(\Omega_m) = 1 -
\sum_{j=1}^{2^{k_m}} \mu(\mathrm{D}_j^{k_m} \setminus
\mathrm{K}_j^{k_m}) \ge 1 - 2^{-k_m} \longrightarrow 1 \qquad
\mbox{as} \qquad m \longrightarrow \infty.$$ Then, the announced
properties follow. Indeed, we have
\begin{itemize}
\item[i)] $\widehat{f}_1, \widehat{f}_2, \ldots$ have pairwise
disjoint supports on $\G$, see the definition of $f_m$.
\item[ii)] If $g \in \mathrm{G}$, we have $$|f_m(g)| \le
|\psi_{k_m}(g)| + |f_m(g) - \psi_{k_m}(g)| < |f_0(g)| +
\varepsilon_m.$$
\item[iii)] If $g \in \Omega_m$, we have $$|f_m(g)| \ge
|\psi_{k_m}(g)| - |f_m(g) - \psi_{k_m}(g)| > |f_0(g)| -
\varepsilon_m.$$
\end{itemize}
This concludes the proof. \fin

\noindent \textbf{Note from the author.} The original purpose of
this paper was to study the sharp Fourier type of $L_p$ spaces
with respect to non-commutative compact groups. This is what I did
in the preliminary version of this paper. Unfortunately, I found
an \emph{stupid} mistake but with serious implications in the
proof. What I present here is the material which remains
\emph{alive} from the previous version.

\bibliographystyle{amsplain}

\begin{thebibliography}{10}
\bibitem {Fo} G.B. Folland, A Course in Abstrac Harmonic
Analysis. Stud. Adv. Math., CRC Press (1995).
\bibitem {GMP} J. Garc\'{\i}a-Cuerva, J.M. Marco and J. Parcet, Sharp
Fourier type and cotype with respect to compact semisimple Lie
groups. Trans. Amer. Math. Soc. \textbf{355} (2003), 3591-3609.
\end{thebibliography}

\end{document}